\documentclass[12pt]{article}
\usepackage{amsxtra,amssymb,amsthm,amsmath,latexsym}

\theoremstyle{plain}
\newtheorem{theorem}{Theorem}

\def\R{{\mathbb R}}

\def\oH{{\overset{\circ}{H}}}
\def\oH1{{\overset{\circ}{H}\kern-.02in{}^1}}

\def\bee{\begin{equation*}}
\def\eee{\end{equation*}}
\def\be{\begin{equation}}
\def\ee{\end{equation}}

\def\R{\mathbb{R}}

\begin{document}

\begin{titlepage}
\title{Inverse problem of potential theory }

\author{A. G. Ramm\\
 Mathematics Department, Kansas State University, \\
 Manhattan, KS 66506-2602, USA\\
ramm@math.ksu.edu }
\date{}
\maketitle\thispagestyle{empty}

\begin{abstract}

P. Novikov in 1938 has proved that if $u_1(x)=u_2(x)$ for $|x|>R$, where $R>0$ is a large number,
$$u_j(x):=\int_{D_j}g_0(x,y)dy,  \quad g_0(x,y):=\frac 1 {4\pi |x-y|},$$
and $D_j\subset \R^3$, $j=1,2,$ $D_j\subset B_R$, are bounded, connected, smooth domains, star-shaped with respect to a common point, then $D_1=D_2$. Here $B_R:=
\{x: |x|\le R\}$.

Our basic results are:  

a)  the removal of the assumption about star-shapeness of $D_j$,

b)  a new approach to the problem,

c)  the construction of counter-examples for 
a similar problem in which $g_0$ is replaced by
$g=\frac {e^{ik|x-y|}}{4\pi |x-y|}$, where $k>0$ is a fixed constant.
 
\end{abstract}
\end{titlepage}

\section{Introduction}\label{S:1}

Suppose there are two bodies $D_j$, $j=1,2,$
uniformly charged with charge density $1$. Let the corresponding potentials be $u_j(x)=
\int_{D_j}\frac {dy}{4\pi |x-y|}$. Assume that $u_1(x)=u_2(x)$ for $|x|>R$, where $R>0$ is a large number. The classical question (inverse problem of potential theory) is: {\em  does this imply that $D_1=D_2$?}

P.Novikov in 1938 (see \cite{N}) has proved a uniqueness theorem for the solution of inverse problem (IP)  of potential theory under a special assumption, see Proposition 1 below.

Let
\be\label{e1}
u(x)=\int_{D}g_0(x,y)dy, \quad g_0(x,y):=\frac 1 {4\pi |x-y|},
\ee
where $D\subset \R^3$ is a bounded, connected, $C^2-$smooth domain. 

We use the following notations: $D_j$, $j=1,2,$ are two different domains $D$,
 $S_j$ is the boundary of $D_j$,  
 $D_j'=\R^3\setminus D_j$,  $S^2$ is the unit sphere in $\R^3$,  $B_R$ is the ball of radius $R$, centered at the origin, $B_R'=\R^3\setminus B_R$, $D_j\subset B_R$, $D_{12}:=D_1\cup D_2$, $D_{12}':=\R^3\setminus D_{12}$, $\mathcal{D}:=D^{12}=D_1\cap D_2$, $\mathcal{D}'=\R^3\setminus \mathcal{D}$,
 $\mathcal{S}$ is the boundary of $\mathcal{D}$, and let
 $u_j(x)=\int_{D_j}g_0(x,y)dy$, $j=1,2.$ 

 P. Novikov has proved the following result, see \cite{N}:
 
 {\bf Proposition 1.} {\em If $u_1(x)=u_2(x)$ for $|x|>R$, then $D_1=D_2$ provided that $D_j$, $j=1,2,$ are star-shaped with respect to a common point. }

In \cite{M} this result is generalized:
the existence of the common point with respect to which $D_1$ and $D_2$ are star-shaped is not
assumed, but $D_j$ are still assumed 
star-shaped.

In \cite{R470}, p.334, a new proof of Proposition 1 was given.  In \cite{R512}--\cite{R190}, see also \cite{R676}, some inverse problems and symmetry problems are studied. 

The goal of this paper is to give a new method for a proof of a generalization of Proposition 1.
In this generalization 

a) The assumptions about star-shapeness of $D_j$, $j=1,2,$  are discarded (see Theorem 1 below); 

b) A new approach to the IP is developed; 

and

c) A similar inverse problem is studied in the case when $g_0(x,y)$ is replaced by the Green's function of the Helmholtz operator, $g(x,y):=\frac{e^{ik|x-y|}}{4\pi |x-y|}$, $k=const>0$ is fixed, and $u_j$ is replaced by 
$$U_j(x)=\int_{D_j} g(x,y)dy.$$ 
The result is formulated in Theorem 2, below.

The idea of our proof does not use the basic idea of \cite{N}, \cite{M},
or \cite{R470}. Our proof is
based on some lemmas.

{\bf Lemma 1.} {\em If $R(\phi)$ is the element of the rotation group in $\R^3$ then $\frac{\partial R(\phi)x}{\partial \phi}|_{\phi=0}=[\alpha, x]$,
where $[x,y]$ is the cross product of two vectors,  $\alpha$ is the unit vector around which the
rotation by the angle $\phi$ takes place, and $x$ is an arbitrary vector.}
 
This lemma is proved in \cite{R470}, p. 416.

{\bf Lemma 2.} {\em The set of restrictions on $S$ of all harmonic in $B_R$ functions is dense
in $L^2(S)$, where $S$ is the boundary of $D$.}

{\em Proof of Lemma 2.} Let us assume the contrary and derive a contradiction. Without loss of generality one may assume $f$ to be real-valued. 
 Suppose $f\not\equiv 0$ is orthogonal in $L^2(S)$ to any harmonic function $h$, that is,
\be\label{e2}
 \int_Sfhds=0
 \ee
 for all harmonic functions in $B_R$.
 
Define $$v(x):=\int_S g_0(x,s)fds.$$
 Assumption  
\eqref{e2} implies $v(x)=0$ in $\R^3$. Indeed,
there exists a unique solution to the problem
$$\Delta h=0\quad in \quad D, \quad h|_S=f.$$ 
For this $h$ equation \eqref{e2} implies that $f=0$, so $v=0$ in $\R^3$.  \hfill$\Box$

{\bf Remark 1.} The proof of Lemma 2 is valid for closed surfaces $S$ which are not necessarily connected. For example, it is valid for $S$ which is a union of two surfaces.

It is known (see, for example, \cite{R470}) that
\be\label{e3}
g_0(x,y)=g_0(|x|)\sum_{\ell \ge 0}\frac {|y|^\ell}{|x|^\ell}\overline{Y_\ell (y^0)}Y_\ell(x^0), \quad |x|>|y|,
\ee
where $x^0:=x/|x|$, $Y_\ell$ are the  spherical harmonics, normalized in $L^2(S^2)$, and $|y|^\ell Y_\ell(y^0)$ are harmonic functions.  The set of these functions for all $\ell\ge 0$ is dense in the set of all harmonic functions in $B_R$. 

{\bf Lemma 3.} {\em If $S$ is a smooth closed surface and $[s,N]=0$ on $S$, then $S$ is a sphere.}

{\em Proof of Lemma 3.} Let $s=s(p,q)$ be a parmetric equations of $S$. Then $N$ is proportional to $[s_p,s_q]$,
where $s_p$ is the partial derivative $\frac {\partial s}{\partial p}$. If $[s,N]=0$, then $[s,[s_p,s_q]]=0$. Thus,
$s_p s\cdot s_q- s_q s\cdot s_p=0$, where $s\cdot s_q$
is the dot product of two vectors. Since $S$ is smooth,
vectors $s_p$ and $s_q$ are linearly independent on $S$.
Therefore, $\frac {\partial s \cdot s}{\partial p}=0$,
and   $\frac {\partial s \cdot s}{\partial q}=0$. Consequently, $s\cdot s=const$. This means that $S$ is a sphere.\\
 Lemma 3 is proved.\hfill$\Box$

Lemma 3 is Lemma 11.2.2 in \cite{R470}, see also Theorem 2 in \cite{R629} and \cite{R666}. Its short proof
is included for convenience of the reader.

It follows from our proof that if $S$ has finitely many
points of non-smoothness, then the parts of $S$, joining these points, are spherical segments.

  Our new results  are the following theorems.
\begin{theorem}\label{T:1}
If $u_1(x)=u_2(x)$ for $|x|>R$, then $D_1=D_2$.
 \end{theorem}

\begin{theorem}\label{T:2}
There may exist, in general, countably many different 
$D_j$ such that the corresponding potentials $U_j$ are equal in $B_R'$ for sufficiently large $R>0$.
 \end{theorem}

{\bf Remark 2.} If  
$$V_j:=\int_{S_j} g(x,t)dt,\qquad j=1,2, $$
then there exist $S_1\neq S_2$ for which $V_1(x)=V_2(x)=0$ $\forall x\in B_R'$.

An example can be constructed similarly to the one given in the proof of Theorem 2.
 
In Section 2 proofs are given.

\section{Proofs}\label{S:2}

\begin{proof}[Proof of Theorem 1.]  If $u_1=u_2$ 
for all $x\in B_R'$ then it follows from the asymptotic of $u_j$ as $|x|\to \infty$ that $|D_1|=|D_2|$. 
Thus, the case $D_1\subset D_2$ is not possible  if $u_1=u_2$ for all $x\in B_R'$. 

  The functions $u_j$ are harmonic functions in $D_j'$, that is, $\Delta u_j=0$ in $D_j'$.
If $u_1=u_2$ in $B_R'$ and $D_1\neq D_2$, then
$u_1=u_2$ in $D_{12}'$
 by the unique continuation property for harmonic functions.

Let $w:=u_1-u_2$. Then $w=0$ in $D_{12}'$,
\be\label{e4}
\Delta w=\chi_2-\chi_1,
\ee
where $\chi_j$ is
the characteristic function of $D_j$. 

Let 
$h$ be an arbitrary harmonic function in $B_R$. Then
$$\int_{D_{2}}h(x)dx-\int_{D_{1}}h(x)dx=0,$$ 
as one gets multiplying  
\eqref{e4} by $h$, integrating by parts and taking into account that $w$ vanishes outside $D_{12}$.

If $h(x)$ is harmonic, so is $h(R(\phi)x)$. Thus,
\be\label{e5}
\int_{D_{2}}h(R(\phi)x)dx-\int_{D_{1}}h(R(\phi)x)dx=0. 
\ee
Differentiate \eqref{e5} with respect to
 $\phi$ and then set $\phi=0$. Using 
 Lemma 1, one gets:
\be\label{e6}
\int_{D_{2}} \nabla h \cdot [\alpha, x]dx-
\int_{D_{1}} \nabla h \cdot [\alpha, x]dx=0,
\ee
where $\alpha\in S^2$ is arbitrary, and $h$ is an arbitrary harmonic function in $B_R$. Since \newline
$ \nabla h \cdot [\alpha, x]=\nabla \cdot (h[\alpha, x])$, it follows from \eqref{e6} and the divergence theorem that
\be\label{e7}
\int_{S_{2}}Nh[\alpha,s]ds-\int_{S_{1}}Nh[\alpha,s]ds=0, 
\ee
for all $\alpha \in S^2$ and all harmonic
 $h$  in $B_R$. Here $N$ is the unit normal to the boundary pointing out of $D_j$.

  If $D_1=D_2$ then \eqref{e7} is an identity. Suppose that $D_1\neq D_2$.
Since $\alpha$ is arbitrary, it follows from \eqref{e7} that
\be\label{e7'}
\int_{S_{2}}h[N,s]ds-\int_{S_{1}}h[N,s]ds=0,  
\ee
for all harmonic in $B_R$ functions $h$.

By Lemma 2 and Remark 1 it follows from \eqref{e7'} that 
$[N,s]=0$ on $S_2$ and on $S_1$. 

By Lemma 3 it follows that $S_1$ and $S_2$ are spheres, so $D_1$ and $D_2$ are balls. These balls must be of the same radius, as was mentioned earlier. 

Now we have a contradiction unless $D_1=D_2$, because two uniformly charged balls  with the same total charge but with different centers cannot have the same potential in $B_R'$. This follows from the explicit formula for their potentials. Theorem 1 is proved. 
\end{proof} 
  
\begin{proof}[Proof of Theorem 2.]
Let $D_j=B_{a_j}$, where $a_j>0$ are some numbers which
are chosen below. Then $U_j(x):=\int_{|y|\le a_j}g(x,y)dy=0$ in the region $B_{a_j}'$ if and only if 
$$\int_0^{a_j}r^2j_0(kr)dr=0,$$
where $j_0(r)$ is the spherical Bessel function,
$$j_0(kr):=\Big(\frac {\pi}{2kr}\Big)^{1/2}J_{1/2}(kr)=\frac{\sin(kr)}{kr}.$$  
This follows from the formula $U_j(x)=\int_{B_{a_j}}g(x,y)dy$, and from the known formula for g(x,y) (see, for example, \cite{R470}):
$$g(x,y)=\sum_{\ell\ge 0}
ik j_\ell(k|y|)h_\ell(k|x|) \overline{Y_\ell(y^0)}Y_\ell(x^0), \quad |y|<|x|.$$
 Here
$Y_\ell$ are the normalized spherical harmonics (see \cite{R470}, p. 261), $j_\ell$ and $h_\ell$ are the spherical Bessel and Hankel functions (see \cite{R470}, p. 262), and the known formula 
 $$\int_{S^2}Y_\ell(y^0)dy^0=0, \qquad \ell>0$$
 was used. 

One has
$$\int_0^{a_j}r^2j_0(kr)dr=\frac{\sin(ka_j)}{k^2}- \frac{a_j\cos(ka_j)}{k^2}=0$$ if and only if 
$$\tan(ka_j)=a_j.$$
 This equation has countably
many positive solutions. To each of these solutions there corresponds a ball $B_{a_j}$ such
that $U_j=0$ in  $B_{a_j}'$. Thus, there are many different balls for which $U_j$ are the same  in $B_R'$, namely $U_j=0$ in $B_R'$ for $R>a_j$. 
Theorem 2 is proved.
\end{proof}


\end{document}